\documentclass{amsart}
\usepackage{amsmath,amsthm,bm,mathrsfs,setspace}
\usepackage[author-year]{amsrefs}

\author{Davar Khoshnevisan}
\title[Schoenberg's Theorem]{Schoenberg's Theorem\\
	Via the Law of Large Numbers\\
	\textbf{(Not for Publication)}}\thanks{Research supported in part by a
	grant from NSF}
\address{144 S 1500 E, Department of Mathematics,
	The University of Utah, Salt Lake City UT 84112--0090}
\email{davar@math.utah.edu}
\urladdr{http://www.math.utah.edu/~{}davar}
\date{April 29, 2005}
\numberwithin{equation}{section}
\begin{document}
\onehalfspacing

\theoremstyle{remark} \newtheorem*{remark}{Remark}
\theoremstyle{theorem} \newtheorem{theorem}{Theorem}
\theoremstyle{theorem} \newtheorem{lemma}{Lemma}

\theoremstyle{theorem} \newtheorem*{schoenberg}{Schoenberg's Theorem}

\newcommand{\R}{\mathbf{R}}

\begin{abstract}
	A classical theorem of S. Bochner states that a function
	$f:\R^n\to\mathbf{C}$ is the Fourier transform of
	a finite Borel measure if and only if $f$ is positive
	definite. In 1938, I. Schoenberg found a beautiful
	complement to Bochner's theorem.
	We present a non-technical derivation of
	of Schoenberg's theorem that relies chiefly on the
	de Finetti theorem and the
	law of large numbers of classical probability theory.
\end{abstract}
\maketitle

\section{Introduction}
A real-valued function $g$ of $n$ vectors is 
said to be \emph{positive semi-definite}
(sometimes, \emph{positive definite})
if  $\sum_{i=1}^k \sum_{j=1}^k g(x_i-x_j) c_i
\overline{c}_j\ge 0$
for all $n$-vectors $x_1,\ldots,x_k$ and 
all complex numbers $c_1,\ldots,c_k$.

A classical theorem of S.  Bochner
\ycite{bochner}*{Theorem 3.2.3, p.\ 58} asserts that
positive semi-definite functions are precisely
those that are Fourier transforms of 
finite measures. Let $\|\cdot\|_n$ denote
the usual Euclidean norm in $n$ dimensions.
That is, $\|x\|_n=(x_1^2+\cdots+x_n^2)^{1/2}$
for all $x\in\R^n$. Then,
the goal of this note is to present a
very simple proof of the following well-known theorem of
I.  J.  Schoenberg~\ycite{schoenberg}*{Theorem 2}:

\begin{schoenberg}
	Suppose $f:\R_+\to\R_+$ is continuous.
	Then, the following are equivalent:
	\begin{enumerate}
		\item The function $\R^n\ni x\mapsto
			f(\|x\|_n)$ is
			positive semi-definite.
		\item The function $\R_+\ni t\mapsto f(\sqrt{t})$ is
			the Laplace transform of a finite Borel measure on $\R_+$.
	\end{enumerate}
\end{schoenberg}

Originally, this theorem 
was used to describe isometric embeddings of Hilbert spaces.
Since its discovery, it has also found non-trivial
connections to other diverse areas ranging
from classical, as well as abstract,
harmonic analysis \citelist{\cite{bergressel}\cite{bcr}\cite{kahane}}
to the measure theory of Banach spaces
\citelist{\cite{bdk1}\cite{bdk2}\cite{bdk3}\cite{CR}\cite{koldobsky2}
\cite{misiewicz1}\cite{misiewicz2}\cite{koldobsky}\cite{koldobskylonke}}, 
function theory \citelist{\cite{ressel:74}\cite{df1}\cite{df2}} and
to the foundations of statistics via de Finetti-type theorems
\citelist{\cite{freedman}\cite{ressel:85}\cite{df1}\cite{df2}}. For other relations,
in particular, to statistical mechanics, see the detailed
historical section of \ocite{df2}.

Although Schoenberg's original proof is not too difficult
to follow, it is somewhat technical.  P.  Ressel
\ycite{ressel:76} has devised a simpler proof
which rests on a
characterization of Laplace transforms
\cite{ressel:74}*{Satz 1} that is similar to
Schoenberg's theorem. We are aware also
of  another simple proof,
due to J.  Bretagnolle, D.  Dacuhna--Castelle,
and J.-L.  Krivine~\citelist{\ycite{bdk1}\ycite{bdk2}\ycite{bdk3}}.
Their proof is similar to the one presented here, but is
slightly more technical.

The present article aims to
describe a self-contained, elementary, and brief derivation
of Schoenberg's theorem. Our proof assumes only
a brief acquaintance with real analysis and
measure-theoretic probability theory.
This proof is quite robust and can be used to produce
more general results;
all one needs is a more general setting in which a
basic form of the de Finetti theorem
and the law of large numbers hold.

Since writing the first draft of this paper,
we have found out about
the work of D. Kelker~\ycite{kelker}*{Theorem 10}.
Kelker's proof is essentially the same as ours.
J. Kingman~\ycite{kingman} contains yet another
rediscovery of Kelker's proof.\\

\noindent\textbf{Acknowledgements.}\
Christian Berg brought
to my attention the recent work of \ocite{sterrneman},
and Paul Ressel made an important correction to the
original draft. I am deeply endebted to them both.

\section{The Proof}
All notation and references to probability theory
are standard and can be found in any
standard first-year graduate textbook. 

Without loss of generality, we may suppose that $f(0)=1$.
Then, thanks to  Bochner's theorem, 
Schoenberg's theorem translates
to the equivalence of the following two assertions:
\begin{enumerate}
	\item[($1^\circ$)] For all $n\ge 1$ there exists a
		Borel probability measure $\mu_n$ on $\R^n$
		such that
		\begin{equation}
			f\left( \sqrt{\sum_{i=1}^n x_i^2}\right) =
			\int_{\R^n} e^{ix \cdot y}\, \mu_n(dy)
			\quad{}^\forall x:=(x_1,\ldots,x_n)\in\R^n.
		\end{equation}
	\item[($2^\circ$)] There exists a Borel probability measure
		$\nu$ on $\R_+$ such that
		\begin{equation}
			f( t ) = \int_0^\infty e^{- t ^2 s/2}\, \nu(ds)
			\qquad{}^\forall   t > 0.
		\end{equation}
\end{enumerate}
Therefore, it suffices to prove that ($1^\circ$) and ($2^\circ$)
are equivalent. The assertion,
``($2^\circ$)$\Rightarrow$($1^\circ$)'' follows from a direct computation
because $\|x\|_n \mapsto \exp(-\|x\|_n^2s/2)$ is manifestly a Fourier
transform on $\R^n$. So
we prove only the converse. Henceforth, we assume
that ($1^\circ$) holds.

Our next lemma follows immediately from ($1^\circ$)
and the uniqueness theorem.

\begin{lemma}\label{lem:key}
	The family $\{\mu_n\}_{n=1}^\infty$ 
	is consistent.
\end{lemma}

It might help to recall that ``$\{\mu_n\}_{n=1}^\infty$ is 
\emph{consistent}'' means that
for all $n\ge 1$
and all linear Borel sets $A_1,A_2,\ldots,$
$\mu_n(A_1\times\cdots\times A_n) = \mu_{n+1}
(A_1\times \cdots A_n\times \R)$.

\begin{proof}[Proof of Schoenberg's Theorem]

	In accord with Lemma~\ref{lem:key}
	and the Kolmogorov consistency theorem,
	there exists an exchangeable stochastic process $\{Y_k\}_{k=1}^\infty$,
	on some probability space $(\Omega,\mathscr{F},\mathrm{P})$,
	such that  for all $n\ge 1$ and all Borel sets $A\subset\R^n$,
	\begin{equation}\label{eq:muY}
		\mathrm{P} \{(Y_1,\ldots,Y_n)\in A\}=\mu_n(A).
	\end{equation}

	Choose and fix some $ t > 0$,
	and introduce a sequence $\{X_i\}_{i=1}^\infty$
	of independent
	random variables such that
	every $X_i$ has the normal distribution with mean
	$0$ and variance $ t ^2$. We can assume, without loss
	of generality, that the $X_i$'s are defined on the same
	probability space $(\Omega,\mathscr{F},\mathrm{P})$.
	We first apply ($1^\circ$) with $x:=n^{-1/2}(X_1,\ldots, X_n)$, and
	then take expectations, to deduce that for all $n\ge 1$,
	\begin{equation}\begin{split}
		\mathrm{E} \left[ f\left( 
			\sqrt{\frac1n \sum_{i=1}^n X_i^2} \right)\right]
			&= \int_{\R^n} \exp\left( -\frac{ t ^2
			\| y \|_n^2}{2n}\right)\, \mu_n(d y )\\
		&=\mathrm{E}
		\left[  \exp\left( -\frac{ t ^2}{2n} \sum_{i=1}^n Y_i^2
		\right) \right].\label{eq:key}
	\end{split}\end{equation}
	See \eqref{eq:muY} for the last identity.
	Now let $n\to\infty$. The simplest form
	of the law of large numbers dictates that
	$\sum_{i=1}^n X_i^2/n\to  \mathrm{Var} X_1 =t ^2$ in probability.
	Therefore,
	the left-hand side of \eqref{eq:key}
	converges to $f( t )$ by the dominated convergence theorem.
	
	By the de Finetti theorem,
	the $Y_i$'s are conditionally i.i.d.\ given 
	the exchangeable $\sigma$-algebra generated
	by the $Y_i$'s.
	Thanks to the Kolmogorov strong law of large numbers,
	and by the Fubini--theorem, 
	$L :=\lim_{n\to\infty} \frac1n \sum_{i=1}^n Y_i^2$
	exists a.s.
	Moreover, the event 
	$\{ L <\infty\}$ agrees upto null sets
	with $\{ \mathrm{E}[Y_1^2\,|\,\mathscr{E}] <\infty\}$,
	where $\mathscr{E}$ denotes the exchangeable $\sigma$-algebra
	of $\{Y_i\}_{i=1}^\infty$.
	By the dominated convergence theorem,
	the right-hand side of \eqref{eq:key} converges
	to $\mathrm{E}[\exp(- t ^2L/2); L<\infty]$.

	We have proved that $f(t)=\mathrm{E}[\exp(-t^2L/2);
	L<\infty]$
	for a possibly-degenerate non-negative random variable
	$L$. Set $t=0$ to find that $L$ is a proper random
	variable; i.e., $1=f(0)=
	\mathrm{P}\{L<\infty\}$.
	Therefore, ($2^\circ$) follows with $\nu$ denoting the
	distribution of $L$.
\end{proof}


\end{document}